\documentclass[a4paper, 11pt]{article}

\usepackage[applemac]{inputenc}
\usepackage{amssymb} 
\usepackage{graphicx}

\def\bR{\mathbb{R}}	
\def\bE{\mathbb{E}}	
\def\I{\mathbb I}	

\newtheorem{theorem}{Theorem}
\newtheorem{proposition}{Proposition}

\newtheorem{lemma}{Lemma}

\newcommand{\EndProof}{{\hfill$\square$\medskip}}
\newenvironment{proof}{\smallskip\noindent
		{\it Proof.} }{\EndProof}

\newcommand{\f}{{\varphi}}

\def\xut{x^{(t)}}

\def\BX{\mathbf{X}}
\def\BY{\mathbf{Y}}
\def\BZ{\mathbf{Z}}

\newcommand{\CH}{{\mathcal H}}

\newcommand{\CK}{{\mathcal K}}

\newcommand{\CN}{{\mathcal N}}
\newcommand{\CO}{{\mathcal O}}

\def\MH{Metropolis-Hastings}

\title{\bf Selection of a MCMC 
simulation strategy\\
via an entropy convergence criterion}

\author{Didier {\sc Chauveau}${}^1$
\and Pierre {\sc Vandekerkhove}${}^2$}

\date{${}^1$ Université d'Orléans \& CNRS,
${}^2$ Université de Marne-la-Vallée \& CNRS\\
\mbox{}\\
May 10th, 2006}

\begin{document}
\maketitle

{\small \noindent {\bf Abstract.}
In MCMC methods, such as the Metropolis-Hastings (MH) algorithm, the
Gibbs sampler, or recent adaptive methods, many different strategies
can be proposed, often associated in practice to unknown rates of
convergence.  In this paper we propose a simulation-based methodology
to compare these rates of convergence, grounded on an entropy
criterion computed from parallel (i.i.d.) simulated Markov chains
coming from each candidate strategy.  Our criterion determines on the
very first iterations the best strategy among the candidates.
Theoretically, we give for the MH algorithm general conditions under
which its successive densities satisfy adequate smoothness and tail
properties, so that this entropy criterion can be estimated
consistently using kernel density estimate and Monte Carlo
integration.  Simulated examples are provided to illustrate this
convergence criterion.}





\section{Introduction}
\label{sec:int}

A Markov Chain Monte Carlo (MCMC) method generates an ergodic Markov
chain $\xut$ for which the stationary distribution is a given probability
density function (pdf)~$f$ over a state space $\Omega\subseteq\bR^s$.
In situations where
direct simulation from~$f$ is not tractable, or where integrals like
$\bE_{f}[h] = \int h(x)f(x)\,dx$
are not available in closed form, MCMC method is
appropriate since, for $T$ large enough, $x^{(T)}$ is approximately
$f$ distributed, and $\bE_{f}[h]$ can be approximated by ergodic
averages from the chain. A major context is Bayesian inference,
where~$f$ is a posterior distribution, usually known only up to a
multiplicative normalization constant.

The Metropolis-Hastings (MH) algorithm (Hastings~\cite{Has})
is one of the most popular
algorithm used in MCMC methods.
Another commonly used MCMC methods is the Gibbs sampler (first
introduced by Geman and Geman~\cite{GemGem}; see also Gelfand and
Smith~\cite{GelSmi}).  An account of definitions and convergence
properties of Gibbs and MH algorithms can be found, e.g., in
Gilks {\it et al.} \cite{GRS}.

In this paper, the theoretical developments will be 
focused on the MH algorithm, since the generic form of its kernel
allows for a general study, as indicated below. However, our proposed methodology
can be applied empirically to any MCMC algorithm (e.g., to the Gibbs sampler). 
It can also be applied to compare the various recent
adaptive methods, which is an area of current and growing
research in MCMC.

For the MH algorithm,
the ``target'' pdf~$f$ needs to be known only up to a (normalizing)
multiplicative constant.
Each step is based on the generation of the proposed next move $y$
from a general conditional density $q(y|x)$, called the {\it instrumental
distribution} or {\it proposal density} (hence a practical requirement is that
simulations from $q$ should be done easily).
For a starting value $x^{(0)}\sim p^{0}$, the $n$-th step
$x^{(n)}\to x^{(n+1)}$ of the algorithm is as follows:

\begin{enumerate}
{\tt
\item[1.] generate $y \sim q(\cdot|x^{(n)})$
\item[2.] compute $\alpha(x^{(n)},y) = \displaystyle
\min \left\{1 , {f(y) q(x^{(n)}|y)\over f(x^{(n)}) q(y|x^{(n)})}
\right\}$
\item[3.] $\displaystyle
\mbox{take }x^{(n+1)} = \left\{
\begin{array}{ll}
y & \mbox{with probability } \alpha(x^{(n)},y),\\ x^{(n)} &
\mbox{with probability } 1-\alpha(x^{(n)},y). \end{array}
\right.$
}
\end{enumerate}

Two well-known MH strategies are (i) the (so-called)
{\it Independence Sampler} (IS), i.e.\ the MH
algorithm with proposal distribution $q(y|x) = q(y)$ independent
of the current position, and (ii) the Random Walk MH algorithm (RWMH),
for which the proposal is a random perturbation $u$ of the current position,
$y = x^{(n)} + u$. The usual choice for the latter is a
gaussian perturbation with a fixed variance matrix
(e.g., in the one-dimensional case, $q(y|x)$ is the pdf of
$\CN(x,\sigma^2)$ where $\sigma^2$ is the scaling parameter of the
perturbation, that has to be tuned).

Ergodicity and convergence properties of the MH algorithm have been
intensively studied in the literature, and conditions have been
given for its geometric convergence (see, e.g., Mengersen and
Tweedie~\cite{MenTwe}, or Roberts and Tweedie~\cite{RobTwe}).
In particular,
Mengersen and Tweedie proved geometric convergence in total variation
norm of the IS, under the condition
$q(y) \ge a f(y)$ for some $a>0$. The associated geometric rate is
$(1-a)^n$, not surprisingly pointing out the
link between the convergence rate and the proximity of~$q$ to~$f$.

To actually implement the MH algorithm, a virtually unlimited number
of choices for the instrumental distribution can be made, with the
goal of improving mixing and convergence properties of the resulting
Markov chain.  If one wants to use the IS strategy, selection of a
reasonably ``good'' proposal density can be done using several
procedures, among which: numerical analysis or a priori knowledge
about the target to approximate the shape of~$f$ (modes,\ldots);
preliminary MCMC experiment, or adaptive methods to dynamically build a
proposal density on the basis of the chain(s)
history (Gelfand and Sahu~\cite{GelSah},
Gilks {\it et al.}~\cite{GRG}, \cite{GiRoSa},
Chauveau and Vandekerkhove~\cite{CV1},
Haario {\it et al.}~\cite{HST}).
If one wants to use the RWMH strategy, ``good'' scaling constants
must be found, since the mixing depends
dramatically on the variance matrix of the perturbation
(see, e.g., Roberts and Rosenthal~\cite{RoRo}).
However, these various choices are associated in general to unknown
rates of convergence, because of the complexity of the kernel, and of
the associated theoretical computations of bounds.

The Gibbs sampler (Geman and Geman~\cite{GemGem})
is defined in a multidimensional setup ($s>1$). 
It consists in simulating a Markov chain 
$x^{(n)}=\bigl(x_1^{(n)},\ldots,x_s^{(n)}\bigr)$ by simulating each 
(not necessarily scalar) coordinate
 according to a decomposition of $f$ in a set 
of its full conditional distributions. In the case of a decomposition 
in $s$ scalar coordinates, the $n$th step of the Gibbs sampler is:
\begin{enumerate}
{\tt
\item[1.]  $x_1^{(n+1)} \sim
           f_1\left( x_1 | x_2^{(n)},\ldots, x_s^{(n)} \right)$
\item[2.]  $x_2^{(n+1)} \sim
           f_2\left(
x_2 | x_1^{(n+1)}, x_3^{(n)},\ldots, x_s^{(n)} \right)$
\item[ ]  \ldots
\item[s.]  $x_s^{(n+1)} \sim
           f_s\left(x_s | x_1^{(n+1)},\ldots, x_{s-1}^{(n+1)} \right)$.
}
\end{enumerate}
Their exists formally many possible decomposition of $f$ in a
set of full conditionals, each of which resulting in a different Gibbs
sampler. In addition, data augmentation schemes 
(see Tanner and Wong~\cite{TW87}) may be used to sample 
from~$f$, which gives even more possibilities,
resulting here also in several simulation
strategies that lead to generally unknown rates of convergence.
Hence our entropy criterion may be used also in this setup to
compare different Gibbs samplers, or to compare Gibbs samplers against MH 
algorithms or other strategies for the same target~$f$.

The motivation of
this paper is thus to propose a method to compare the rates of
convergence of several candidate simulation algorithms  (designed for the
same target~$f$), solely on the basis of the
simulated output from each Markov chain.
Note that the question of selecting the best MH strategy among a
family of proposal densities is the subject of recent developments
(see, e.g., G\aa semyr~\cite{GasJor} for a heuristical solution using
adaptation, Mira~\cite{Mi01} for an ordering of MCMC algorithms
based on their asymptotic precisions, or 
Rigat~\cite{Ri06}).

We suggest the use of an entropy
criterion between the pdf of each algorithm
at time~$n$ and the target density~$f$.  The
computation of an estimate of this criterion requires the simulation,
for a short duration~$n_0$,
of~$N$ parallel (i.i.d.) chains coming from each strategy.
This can be seen as a pre-run, to determine the best algorithm before
running it for the long duration required by the MCMC methodology.

More precisely, assume we have two simulation
strategies, generating two MCMC algorithms with densities denoted by
$p_1^n$ and $p_2^n$ at time (iteration) $n$.
For the comparison, both algorithms are started with the same
initial distribution, i.e. $p_1^0 = p_2^0$.
Define the relative entropy of a probability density~$p$ by
\begin{equation}
\label{eq:hp}
\CH(p)=\int p(x)\log p(x)\, dx.
\end{equation}
A natural measure
of the algorithm's quality is the evolution in time ($n$) of the
Kullback-Leibler ``divergence'' between $p^n_i$, $i=1,2,$ and $f$, given by
$$
{\cal K}(p^n_i,f) = \int  \log\left( {p^n_i(x)\over f(x)} \right)
             p^n_i(x)\,dx
             = \CH(p^n_i) - \bE_{p^n_i}[\log f].
$$
The behavior of the application $n\to{\cal K}(p^n_i,f)$ will be
detailed in section~\ref{sec:Kull}.

When $f$ is analytically known, an a.s.\ consistent estimation of
$\bE_{p^n_i}[\log f]$ is obtained easily by Monte-Carlo integration
using $N$ i.i.d.\ realisations from the algorithm at time $n$.
Unfortunately in the MCMC setup, $f$ is usually the posterior density
of a Bayesian model, so that $f(\cdot)=C \varphi(\cdot)$ where the
normalization constant~$C$ cannot be computed, hence this direct
estimation cannot be done.  However, if we want to compare two
strategies, knowing $C$ is not needed to estimate the {\it difference}
of the divergences with respect to~$f$. Define
\begin{eqnarray}
D(p_1^n,p_2^n,f) &=& \CK(p_1^n,f) - \CK(p_2^n,f)\nonumber\\
&=&\CH(p_{1}^n) - \CH(p_{2}^n)
+ \bE_{p_2^n}[\log\varphi] - \bE_{p_1^n}[\log \varphi].
\label{eq:Entdiff}
\end{eqnarray}
The Kullback criterion is the only divergence insuring this property,
hence it motivates our choice for applying it.  One may think of other
distances, such as $L^1$ or $L^2$, but estimating such distances
requires regularity conditions similar to ours (see, e.g.,
Devroye~\cite{Dev83}).  In addition, using other divergences would
require an estimation of $C$ by other techniques, which is typically
not feasible in actual MCMC situations.  Note also that the Kullback
divergence is currently used as a criterion in other simulation
approaches (see Douc~{\it et al.\/}~\cite{DGMR05}).

We propose to use the
$N$  i.i.d.\ simulations from each of the two strategies at time~$n$,
$$
     (X_1^{i,n},\ldots,X_N^{i,n}) \quad\mbox{i.i.d.}\sim p_i^n, i=1,2.
$$
These simulations are first used to
estimate $\bE_{p_i^n}[\log \varphi]$ via Monte-Carlo integration.
Denote these estimates by
\begin{equation}
p_i^n(\overline{\log\varphi})_N = {1\over N}\sum_{j=1}^N \log
\varphi(X_j^{i,n})
\stackrel{a.s.}{\longrightarrow} \bE_{p_i^n}[\log\varphi],
\label{eq:MClogphi}
\end{equation}
where the convergence comes from the strong law of large numbers.

The remaining problem is then the estimation of the entropies
$\CH(p^n_i)$, $i=1,2$.
One classical approach is to
build a nonparametric kernel density estimate of $p_i^n$, and to
compute the Monte-Carlo integration of this estimate.
Techniques based on this approach have been suggested by Ahmad and
Lin~\cite{Ahm1},
and studied by several authors under different assumptions
(see, e.g., Ahmad and Lin~\cite{Ahm2},
Eggermont and LaRiccia~\cite{Egger}, Mokkadem~\cite{Mokka}).
An interesting point is that in our setup, we can ``recycle''
the simulations already used to compute
the $p_i^n(\overline{\log\varphi})_N$'s. We denote in the sequel
our entropy estimate of $\CH(p_i^n)$ by $\CH_N(p_i^n)$, which
will be defined and studied in Section~\ref{sec:entropyest}.
We define accordingly
$$
D_N(p_1^n,p_2^n,f) = \CH_N(p_{1}^n) - \CH_N(p_{2}^n)
+ p_2^n(\overline{\log\varphi})_N - p_1^n(\overline{\log\varphi})_N.
$$

Our methodology can be applied in actual situations in the following
way: assume we have $k$ possible simulation strategies
$s_1,\ldots,s_k$ to sample
from~$f$, 
resulting in~$k$ successive densities $p_i^n$, $i=1,\ldots,k$, $n\ge0$.
Let $p_i^0=p^0$ be the common initial distribution for the~$k$ algorithms.
The determination of the best algorithm among the~$k$ candidates
can be done using the steps below:
\begin{enumerate}
     \item  select the best strategy $s_b$ between $s_1$ and $s_2$, on
     the basis of the sign of (the plot of)
     $n\mapsto D_N(p_1^n,p_2^n,f)$, for $n=1,\ldots,n_0$, where $n_0$ is the
     simulation duration;

     \item  store the sequence of estimates
     $\{\CH_N(p_b^n)-p_b^n(\overline{\log\varphi})_N\}_{1\le n\le n_0}$;

     \item  for $i=3,\ldots,k$,

     \begin{enumerate}
	\item select the best strategy between
	$s_b$ and $s_i$, as in step 1. Notice that the computation of
	$D_N(p_b^n,p_i^n,f)$ just requires now that of
	$\{\CH_N(p_i^n)-p_i^n(\overline{\log\varphi})_N\}_{1\le n\le n_0}$;
	\item update $b$ and the sequence
	$\{\CH_N(p_b^n)-p_b^n(\overline{\log\varphi})_N\}_{1\le n\le n_0}$.
     \end{enumerate}
\end{enumerate}

In practice, $n_0$ can be chosen small, since the best strategy is
usually determined during the very first iterations. For the 
one or two dimensional  examples simulated
in Section~\ref{sec:simul}, we have observed that values of
$n_0$ between $10$ and~$30$ is often sufficient.
The point is that the difference between the entropy
contraction rate of each 
strategy is obvious at the very first iterations.

Storing at each step the sequence of estimates of the best strategy
$\{\CH_N(p_b^n)-p_b^n(\overline{\log\varphi})_N\}_{1\le n\le n_0}$
clearly saves computing time; the total number of simulations required
is thus $N\times k\times n_0$. Concerning the computer investment,
a C program for doing the parallel
(i.i.d.) simulations together with entropy estimation for a generic MH
algorithm is available (from the first author) as a starting point.

As stated previously, the technical part of the paper focus on the MH 
algorithm.
Section~\ref{sec:Kull} outlines some links between
ergodicity and convergence to zero of ${\cal K}(p^n,f)$ as $n\to\infty$.
In Section~\ref{sec:smooth}, we establish assumptions on the
proposal density~$q$, $f$ and the initial density $p^0$
to insure that, at each time $n$, adequate smoothness
conditions hold for the successive densities $p^n$, $n\ge 0$.
These conditions are stated for the general (multi--dimensional)
case, and detailed precisely in the Appendix for the one-dimensional situation,
for which practical examples are given.
In Section~\ref{sec:entropyest}, these conditions are used to
define an estimate of $\CH(p^n)$ on the basis of the
i.i.d.\ simulation output. Finally, Section~\ref{sec:simul}
illustrates the behavior of our methodology
for synthetic one and two-dimensional examples.

We provide in Section~\ref{sec:entropyest} theoretical conditions
under which our criterion is proved to converge, and check in the
appendix that these conditions are satisfied in some classical simple
situations, to show that it can reasonably be expected to be a good
empirical indicator in general situations for which the technical
conditions are hard to verify.  However, it is important to insist on
the fact that, from the methodological point of view, our comparison
criterion may be applied to far more general MCMC situations than the
MH algorithm.  For example, the homogeneous Markov property of the
simulated processes does not play any role in the convergence of the
entropy estimates of $\CH(p^n_i)$, since these estimates are based on
i.i.d.\ copies at time~$n$.  Hence our methodology may be applied to
compare the different adaptive sampling schemes proposed in recent
literature (see, e.g., Haario {\it et al.}~\cite{HST}, Atchad\'e and
Rosenthal~\cite{AR05}, Pasarica and Gelman~\cite{PG05}).  Indeed, we
empirically used a preliminary version of this criterion to evaluate
the adaptive MCMC method proposed in Chauveau and
Vandekerkhove~\cite{CV1}, and we apply it successfully in
Section~\ref{sec:simul} to another adaptive MH algorithm.

\section{Kullback divergence to stationarity}
\label{sec:Kull}

In this section we show a property of the evolution in time of
the Kullback-Leibler divergence between the distibutions $p^n$ of the
MH algorithm and the target distribution $f$. It has been proved
(see, e.g., Miclo ~\cite{Miclo}) that for countable
discrete Markov chains, the Kullback-Leibler divergence between the
measure at time~$n$ and the stationary measure with density~$f$
(also denoted~$f$) decreases with
time, i.e. that ${\cal K}(mP,f)\leq {\cal K}(m,f)$, where $mP$ is
the transportation of a measure~$m$ with the Markov kernel~$P$,
defined by $mP(\cdot) = \int m(dx) P(x,\cdot)$.

We denote in the sequel the supremum norm of a real-valued
function $\f$
by
\begin{equation}
          || \f||_{\infty} := \sup_{x\in\Omega}| \f(x) |.
          \label{eq:normsup}
\end{equation}
We first recall a result due to Holden~\cite{Hol} assessing the
geometric convergence of the MH algorithm under a uniform
minoration condition:

If there exists $a\in(0,1)$ such that
$q(y|x)\ge a f(y)$ for all $x,y \in\Omega$, then:
\begin{equation}
         \forall y\in\Omega,\quad
             \left\vert \frac{p^{n}(y)}{f(y)} - 1\right\vert \le
(1-a)^n \left|\left| \frac{p^0}{f} - 1 \right|\right|_{\infty} .
\label{eq:vn}
\end{equation}

We use this result to show that the Kullback-Leibler divergence
between $p^n$ and $f$ decreases geometrically fast in this case:

\begin{proposition}
\label{taux}
If the proposal density of the \MH~algorithm satifies
$q(y|x)\ge a f(y)$, for all $x,y\in\Omega$, and $a\in(0,1)$, then
\begin{eqnarray}
{\cal K}(p^{n},f) \le \kappa \rho^n (1 + \kappa \rho^n) ,
\end{eqnarray}
where $\kappa=||p^{0}/f - 1||_{\infty}>0$, and $\rho=(1-a)$.
\end{proposition}

\begin{proof}
Using equation~\ref{eq:vn}, we have:
\begin{eqnarray*}
{\cal K}(p^{n},f)&=&
\int \log\left(\frac{p^{n}(y)}{f(y)}\right)p^{n}(y)\,dy\\
&\leq& \int \log\left(
\left| \frac{p^{n}(y)}{f(y)}-1\right|+1\right)
\left( \left|\frac{p^{n}(y)}{f(y)}-1\right|+1\right)f(y)\,dy\\
&\le&  \log(\kappa \rho^n + 1) (\kappa \rho^n + 1)
            \le \kappa \rho^n (\kappa \rho^n + 1).
\end{eqnarray*}
\end{proof}

More generally, the question of the convergence in entropy of a Markov
process is an active field of current research; see, e.g.,
Ball {\it et al.}~\cite{BBN}, Del Moral {\it et al.}~\cite{Delmo1}, and
Chauveau and Vandekerkhove~\cite{CV2}.

\section{Smoothness of MCMC algorithms densities}
\label{sec:smooth}

For estimating the entropy of a MCMC algorithm successive
densities~$p^n$, $n\ge0$, we have to check that appropriate
smoothness and tails technical conditions on these successive densities hold.
In our setup, it appears tractable to apply
results on entropy estimation based on a {\it Lipschitz condition}.
Remember that a function $\f: \Omega \to \bR$ is called
$c$-Lipschitz if there
exists a constant $c>0$ such that, for any $y,z \in\Omega$,
$|\f(y) - \f(z)| \le c ||y-z||$.

As stated in the introduction, we will essentially focus on the MH
case, essentially because its kernel is ``generic'',
depending only on $q$ and $f$.
However, there is a major difficulty in this case, coming from the
fact that the MH kernel has a point mass at the current position.  

The difficulty for the Gibbs sampler  is that 
its successive densities  are given by
$$
p^{n+1}(y) = \int p^n(x) g(x,y)\,dx,
$$
where $g$ is the density of the Gibbs kernel, 
\begin{equation}
g(x,y) = f_1(y_1|x_2,\ldots,x_s)
\times f_2(y_2|y_1,x_3,\ldots,x_s)
\times \cdots
\times f_s(y_s|y_1,\ldots,y_{s-1})
\end{equation}
and Lipschitz condition on $p^n$ depends heavily of the decomposition 
of~$f$.  We indeed obtained a Lipschitz conditionfor the first
iterations in the case of a
toy-size Gibbs sampler (s=2), but stating conditions at a reasonably
general level seems not possible.

\subsection{The MH Independence Sampler case}
\label{ss:is}

From the description of the MH algorithm in Section~\ref{sec:int},
we define the off-diagonal transition density of the MH kernel
at step $n$ by:
\begin{eqnarray}\label{taub}
p(x,y)=\left\{
\begin{array}{ll}
q(y|x)\alpha(x,y)&\mbox{if}~x\neq y,\\
0              &\mbox{if}~x=y,
\end{array}\right.
\end{eqnarray}
and set the probability of staying at~$x$,
$$
r(x)=1-\int p(x,y)dy.
$$
The MH kernel can be written as:
\begin{equation}
               \label{eq:MHkernel}
               P(x,dy) = p(x,y)dy + r(x)\delta_x(dy),
\end{equation}
where $\delta_x$ denotes the point mass at $x$.

We focus first on the IS case ($q(y|x) \equiv q(y)$)
since it allows for simpler conditions.
We will see that the minorization condition $q(x)\ge a f(y)$ which
implies geometric convergence of the IS is also needed for our
regularity conditions. One may argue that, in this case, it is also
possible to use an importance sampling scheme (see, e.g.,
Douc~{\it et al.\/}~\cite{DGMR05}). This strategy garanties i.i.d.\
simulated values for~$f$, but requires the normalization of the
estimate (since the normalization constant~$C$ is unknown), which may 
lead to large variance.

Let $p^{0}$ be the density of the initial distribution of the MH
algorithm, which will be assumed to be ``sufficiently smooth'', in a
sense that will be stated later.
We will assume also that the proposal density $q$ and the target
p.d.f.~$f$ are also sufficiently smooth.

From~(\ref{eq:MHkernel}),
the successive densities of the IS are given by the recursive formula
\begin{eqnarray}
p^{n+1}(y) &=& q(y)\int p^n(x)\alpha(x,y)\,dx +
p^n(y) \int q(x)(1-\alpha(y,x))\,dx \label{eq:pn1} \\
&=& q(y) I_{n}(y) + p^n(y)\left(1-I(y)\right), \label{eq:pn}
\end{eqnarray}
where
\begin{eqnarray}
I_n(y) &=& \int p^n(x)\alpha(x,y)\,dx, \quad n\ge 0 \label{eq:In} \\
I(y) &=& \int q(x)\alpha(y,x)\,dx . \label{eq:I}
\end{eqnarray}
For convenience, we introduce the notations
$$
\alpha(x,y) = \phi(x,y) \wedge 1,\quad
\phi(x,y) = {h(x)\over h(y)},\quad
h(x) = {q(x)\over f(x)}.
$$

We consider the first iteration of the algorithm.
From~(\ref{eq:pn}), the regularity properties of the density
$p^{1}$ are related to the
regularity properties of the two parameter-dependent integrals $I_1$
and $I$. Regularity properties of such integrals
are classically handled by the
theorem of continuity under the integral sign (see, e.g.,
Billingsley~\cite{Bil} Theorem~16.8 p.~212). Continuity is
straightforward here:

\begin{lemma}
      \label{lem:cont}
      If $q$ and $f$ are strictly positive and continuous on
      $\Omega\subseteq\bR^s$,
      and $p^0$ is
      continuous, then $p^n$ is continuous on $\Omega$ for $n\ge 1$.
\end{lemma}

\begin{proof}
It suffices to prove continuity for $p^1(y)$ at any~$y_0\in \Omega$.
The integrand of $I_0$, $p^0(x)\alpha(x,y)$,
is continuous in~$y$ at~$y_0$ for any $x\in \Omega$, and
$$
|p^0(x)\alpha(x,y)| \le p^0(x), \quad \forall x,y\in\Omega.
$$
Then $I_0$ is continuous at~$y$ by the Lebesgue's dominated
convergence theorem (since $\Omega$ is a metric space, so that
continuity can be stated in term of limit of sequence).
The same reasonning
applies to $I(y)$ by using $q$ for the dominating function.
\end{proof}

From equation (\ref{eq:pn}), we have directly that
\begin{eqnarray}
|p^{n+1}(y) - p^{n+1}(z)| & \le & ||q||_{\infty}\, |I_n(y) - I_n(z)|
           + ||I_n||_{\infty}\, |q(y) - q(z)|  \nonumber \\
& + & ||1-I||_{\infty} \,|p^n(y) - p^n(z)| \nonumber \\
& + & ||p^n||_{\infty} \,|I(y) - I(z)|, \label{eq:pndecomp}
\end{eqnarray}
so that, to prove recursively that $p^{n+1}$ is Lipschitz, we have
first to prove that $I_n$ and $I$ are both Lipschitz.

\begin{lemma}
      \label{lem:LipIn}
      If $f/q$ is $c_1$-Lipschitz, and
      $\int p^0 h < \infty$, then for all $n\ge1$:
      \begin{enumerate}
      \item[(i)] $\int p^n h < \infty$;
      \item[(ii)] $I_n$ is $(c_1\int p^{n} h)$-Lipschitz.
      \end{enumerate}
\end{lemma}

\begin{proof}
          first we have to check that $\int p^0 h < \infty$ can be iterated.
          This comes directly from the recursive definition~(\ref{eq:pn1})
          (since $0\le r(x)\le 1$):
          \begin{eqnarray*}
	   \int p^1(y)h(y)\,dy & = &
	 \int \left[\int p^0(x)p(x,y)\, dx + p^0(y)r(y) \right]
	 h(y) \,dy\\
	    & \le & \int {q(y)^2\over f(y)}
	    \left[\int p^0(x)\phi(x,y)\, dx\right]
	    \,dy + \int p^0(y) {q(y)\over f(y)} \,dy \\
	    & = & 2 \int p^0(y) h(y) \,dy < \infty .
          \end{eqnarray*}
Hence $\int p^0 h < \infty \Rightarrow \int p^n h < \infty$ for $n\ge1$.
Then, we have
\begin{eqnarray*}
|I_n(y) - I_n(z)|  & \le & \int p^{n}(x)|\alpha(x,y)-\alpha(x,z)|\,dx  \\
           & \le &  \int p^{n}(x)|\phi(x,y)-\phi(x,z)|\,dx  \\
           & \le &  \int p^{n}(x) h(x)
           \left|{f(y)\over q(y)} - {f(z)\over q(z)}\right|\,dx \\
            &\le& \left( c_1\int p^{n} h \right) ||y-z|| .
\end{eqnarray*}
\end{proof}

Note that the hypothesis that $f/q$ is Lipschitz is reasonable
in the IS context. Indeed, one has to choose a proposal density~$q$ with
adequate tails for the MH to be efficient, i.e.\ to converge quickly.
As recalled in the introduction, it has been proved that the IS is
uniformly geometrically ergodic if $q(y)\ge a f(y)$ for some
$a>0$ (Mengersen and Tweedie~\cite{MenTwe}). Actually, these authors
also proved that the IS
is not even geometrically ergodic if this condition is not satisfied.
But satisfying this minoration condition requires $q$ to have tails
heavier than the tails of the target~$f$. Hence, common choices for
implementing the IS make use of heavy-tailed proposal densities (e.g.,
mixtures of multidimensional Student distributions with small degrees
of freedom parameters),
so that $f/q$ is typically a
continuous and positive function which goes to zero when
$||x||\to\infty$. It can then reasonably be assumed to be Lipschitz.
This condition in
lemma~\ref{lem:LipIn} may thus be viewed as
a consequence of the following assumption, which will be used below:

\noindent
{\bf Assumption A:}
{\it $q$ and $f$ are strictly positive and continuous
densities on $\Omega$, and $q$ has heavier tails than $f$, so that
$\lim_{||y||\to\infty}h(y) = +\infty$.
}

We turn now to the second integral $I(y) = \int q(x)\alpha(y,x)\,dx$.
The difficulty here comes from the fact that the integration variable
is now the {\it second} argument of $\alpha(\cdot,\cdot)$.
Hence, applying the majoration used previously gives
$$
|I(y) - I(z)| \le \int q(x)|\phi(y,x) - \phi(z,x)|\,dx
= \int f(x)|h(y) - h(z)|\,dx ,
$$
and since we have made the ``good'' choice for the proposal density
(assumption A), $h = q/f$ is obviously {\it not} Lipschitz.

A direct study of
$\alpha(\cdot,x) = [h(\cdot)/h(x)] \wedge 1$, as it appears in
$I(y)$ (equations~(\ref{eq:pn}) and~(\ref{eq:I}))
is needed here.
Consider a fixed $x\in\Omega$ in the sequel.
Clearly, there exists by~(A) a compact set
$K(x)$ such that for any $y\notin K(x)$, $h(y) \ge h(x)$.
This entails that
$$
\forall y \notin K(x), \quad \alpha(y,x) = 1.
$$
Now, for any $y\in K(x)$, $\alpha(y,x)$ is a continuous function
truncated at one, so that it is uniformly continuous. If we
assume slightly more, i.e.\ that $\alpha(\cdot,x)$
is $c(x)$-Lipschitz, we have proved the following Lemma:
\begin{lemma}
\label{lem:LipI}
If assumption A holds, and if for each $x$ there exists $c(x)<\infty$
such that
\begin{equation}
\forall y,z \in K(x), \quad
|\alpha(y,x) - \alpha(z,x)| \le c(x) ||y-z||,
\label{eq:condLipalpha}
\end{equation}

where $c(x)$ satisfies
\begin{equation}
          \int q(x)c(x)\,dx < \infty,
          \label{eq:condcq}
\end{equation}
then $I$ satisfies the Lipschitz condition:
$$
\forall (y,z) \in\Omega^2,\quad
|I(y) - I(z)| \le \left(\int q(x)c(x)\,dx\right)  ||y-z||.
$$
\end{lemma}
Examples where lemma~\ref{lem:LipI} holds will be given in the
Appendix, for the one-dimensional situation.

\begin{proposition}
\label{prop:pnLip}
If conditions of Lemmas \ref{lem:cont}, \ref{lem:LipIn} and
\ref{lem:LipI} hold, and if
\begin{enumerate}
          \item[(i)] $|| q ||_{\infty} = Q < \infty$ and $q$ is $c_q$-Lipschitz;
          \item[(ii)]  $||p^0||_{\infty} = M < \infty$ and $p^0$ is
$c_0$-Lipschitz;
\end{enumerate}
then the successive densities of the Independance Sampler
satisfy a Lipschitz condition, i.e. for any $n\ge0$, there exists
$k(n)<\infty$ such that
\begin{equation}
\forall (y,z)\in\Omega^2,\quad
|p^n(y) - p^n(z)| \le k(n)\, ||y-z||.
\label{eq:pnLip}
\end{equation}
\end{proposition}

\begin{proof}
Using equation (\ref{eq:pndecomp}), and the fact that
$$
||I_n||_{\infty} \le \int p^n(x)\,dx = 1, \quad
||I ||_{\infty} \le \int q(x)\,dx = 1,
$$
and
\begin{eqnarray*}
||p^n||_{\infty} &\le&
Q ||I_{n-1}||_{\infty}
+ ||p^{n-1}||_{\infty} \, ||1-I(y)||_{\infty} \\
&\le& nQ + M ,
\end{eqnarray*}
we obtain
\begin{eqnarray*}
|p^{n+1}(y) - p^{n+1}(z)| & \le & Q |I_n(y) - I_n(z)|
           +  |q(y) - q(z)|  \\
&  & + |p^n(y) - p^n(z)|
           + (nQ + M) |I(y) - I(z)|.
\end{eqnarray*}
Thus, applying this recursively, (\ref{eq:pnLip}) is satisfied, with
\begin{eqnarray*}
k(n) &=& Q c_1 \int p^n(x) h(x)\,dx + c_q \\
           & &  + ((n-1)Q+M)\int q(x)c(x)\,dx + k(n-1), \quad n\ge 2\\
k(1) &=& Q c_1 \int p^0(x) h(x)\,dx + c_q
             + M\int q(x)c(x)\,dx + c_0 .
\end{eqnarray*}
\end{proof}

\subsection{The general Metropolis-Hastings case}
\label{ss:general}

When the proposal density is of the general form $q(y|x)$ depending on
the current position of the chain, the successive densities of the
MH algorithm are given by
\begin{eqnarray}
               p^{n+1}(y) &=& \int p^n(x)q(y|x)\alpha(x,y)\,dx +
               p^n(y) \int q(x|y)(1-\alpha(y,x))\,dx\nonumber\\
               &=& J_n(y) + p^n(y)\left(1 -
               J(y)\right),\label{eq:pngeneral}
\end{eqnarray}
where
\begin{eqnarray}
J_n(y) &=& \int p^n(x)q(y|x)\alpha(x,y) \,dx, \label{eq:Jn} \\
     J(y) &=& \int q(x|y)\alpha(y,x) \,dx. \label{eq:J}
\end{eqnarray}
In comparison with the IS case, the continuity already requires some
additional conditions. Let $B(y_0,\delta)$ denotes the closed ball
centered at $y_0\in\Omega$, with radius~$\delta$.

\begin{lemma}
\label{lem:contgeneral}
If $q(x|y)$ and $f$ are
           strictly positive and continuous everywhere on both variables, and
           $p^0$ is continuous, and if:
\begin{enumerate}
           \item[(i)] $\sup_{x,y}q(x|y) \le Q < \infty$ ;

           \item[(ii)] for any $y_0\in\Omega$ and some $\delta>0$,
           $\sup_{y\in B(y_0,\delta)} q(x|y)
           \le \varphi_{y_0,\delta}(x)$,
           where $\varphi_{y_0,\delta}$ is integrable;
\end{enumerate}
then $p^n$ is continuous on $\Omega$ for $n\ge 1$.
\end{lemma}

\begin{proof}
As for Lemma~\ref{lem:cont}, it is enough to check the dominating
conditions of, e.g., Billingsley~\cite{Bil}, p.212.
However, for $J$, we need the local condition~(ii) to prove the
continuity of $J(y)$ at any $y_0\in\Omega$.
\end{proof}

Note that the additional condition~(ii) is reasonable. For instance, we
refer to the most-used case of the RWMH with gaussian perturbation of
scale parameter $\sigma^2>0$. In the one-dimensional case,
$q(x|y)$ is the pdf of $\CN(y,\sigma^2)$ evaluated at~$x$, and
one can simply take for condition~(ii)
\begin{eqnarray}
\varphi_{y_0,\delta}(x) &=&
q(x|y_0 - \delta)\I_{ x < y_0 - \delta} +
q(y_0 - \delta|y_0 - \delta)
\I_{[y_0 - \delta, y_0 + \delta]}(x)\nonumber\\
& & + q(x|y_0 + \delta)\I_{x > y_0 + \delta},\label{eq:exphi}
\end{eqnarray}
(i.e.\ the tail of the leftmost gaussian pdf on the left side, the tail
of the rightmost gaussian pdf on the right side, and the value of the
gaussian at the mode inside $[y_0 - \delta, y_0 + \delta]$).

To prove that the successive densities $p^n$ of the general MH
algorithm are Lipschitz, we proceed using conditions at a higher level
than for the IS case, because the successive densities are more
complicated to handle

\begin{proposition}
\label{prop:LipgeneralMH}
If conditions of Lemma \ref{lem:contgeneral} hold, and if
\begin{enumerate}
\item[(i)]  $||p^0||_{\infty} = M < \infty$ and $p^0$ is $c_0$-Lipschitz;

\item[(ii)] $q(\cdot|x)\alpha(x,\cdot)$ is $c_1(x)$-Lipschitz,
with $\int p^n(x) c_1(x)\,dx <\infty$,

\item[(iii)] $J(\cdot)$ is $c_2$-Lipschitz,
\end{enumerate}
then the successive densities of the general MH
satisfy a Lipschitz condition, i.e. for any $n\ge0$, there exists
$\ell(n)<\infty$ such that
\begin{equation}
\forall (y,z)\in\Omega^2,\quad |p^n(y) - p^n(z)| \le \ell(n)\, ||y-z||.
\label{eq:pnLipgeneral}
\end{equation}
\end{proposition}

\begin{proof}
First, it is easy to check
that, similarly to the IS case,
$|| J_n||_{\infty} \le Q$, $||J||_{\infty} \le 1$,
and $||p^n||_{\infty} \le nQ + M$. Then, using the decomposition
\begin{eqnarray*}
|p^{n+1}(y) - p^{n+1}(z)|
&\le& |J_n(y)-J_n(z)| + 2|p^n(y)-p^n(z)| \\
& & + ||p^n||_{\infty}|J(y)-J(z)| ,
\end{eqnarray*}
equation~(\ref{eq:pnLipgeneral})
is clearly a direct consequence of conditions
(ii) and (iii), and the $\ell(n)$'s can be determined recursively
as in the proof of Proposition~\ref{prop:pnLip}.
\end{proof}

Proposition~\ref{prop:LipgeneralMH} may look artificial
since the conditions are clearly ``what is needed''
to insure the Lipschitz property for $p^n$.  However, we show
in the Appendix (section~\ref{ss:ApGeneral})
that these conditions are reasonable, in the sense that they are
satisfied, e.g., in the one-dimensional case for usual RWMH
algorithms with gaussian proposal densities.

\section{Relative entropy estimation}
\label{sec:entropyest}

Let $\BX_N=(X_1,\ldots,X_N)$ be an i.i.d.\ $N$-sample of random
vectors taking values in $\bR^s$, $s\geq 1$, with common probability
density function $p$.  Suppose we want to estimate the relative
entropy of $p$, $\CH(p)$ given by~(\ref{eq:hp}),
assuming that it is well defined and finite.  Various
estimators for $\CH(p)$ based on $\BX_N$
have been proposed and studied in the literature,
mostly for the case $s=1$.  One method to estimate $\CH(p)$
consists in obtaining a suitable density estimate $\hat p_N$ for $p$,
and then susbtituting $p$ by $p_N$ in an entropy-like functional of
$p$.  This approach have been adopted by Dmitriev and Tarasenko
\cite{Dmitri1}\cite{Dmitri2}, Ahmad and Lin \cite{Ahm1}\cite{Ahm2},
Gy\"{o}rfi and Van Der Meulen \cite{Gyor1}\cite{Gyor2}, and Mokkadem
\cite{Mokka} who prove strong consistency of their estimators in
various framework.  More recently Eggermont and LaRiccia~\cite{Egger}
prove, that they get the best asymptotic normality for the Ahmad and
Lin's estimator for $s=1$, this property being lost in higher
dimension.  Another method used to estimate $\CH(p)$ is based on
considering the sum of logarithms of spacings of order statistics.
This approach was considered by Tarasenko \cite{Taras}, and Dudewicz and
Van Der Meulen~\cite{Dude}.

In our case, and due to the rather poor smoothness properties
that can be proved for
the densities $p^n$ we have to consider, we use
the entropy estimate proposed by
Gy\"{o}rfi and Van Der Meulen \cite {Gyor2}, but with smoothness
conditions of Ivanov and Rozhkova~\cite{Iva}: a Lipschitz condition
which appeared tractable in our setup, as shown in
Section~\ref{sec:smooth}.

Following Gy\"{o}rfi and Van Der Meulen \cite {Gyor2}, we
decompose the sample $\BX_N$ into two
subsamples $\BY_N=\left\{Y_i\right\}$ and
$\BZ_N =\left\{Z_i\right\}$, defined by
\begin{eqnarray}
Y_i&=&X_{2i}\quad\mbox{for } i=1,\ldots,[N/2]\label{Y}\\
Z_i&=&X_{2i-1} \quad\mbox{for }i=1,\ldots,[(N+1)/2]\label{Z},
\end{eqnarray}
where $[N]$ denotes the largest integer inferior to $N$.

Let $\hat p_N(x)= \hat p_N(x,\BZ_N)$
be the Parzen-Rosenblatt kernel density
estimate given by
\begin{eqnarray}
\label{parzen}
\hat p_{N}(x)=\frac{1}{h_N^s (N+1)/2}\sum_{i=1}^{[(N+1)/2]}
K_{h_N}\left(\frac{x-Z_i}{h_N}\right),\quad x\in\bR^s,
\end{eqnarray}
where the kernel $K$ is a density and $h_N>0$ with
$\lim_{N\rightarrow \infty} h_N=0$, and
$\lim_{N\rightarrow \infty}N h^s_N=\infty$.
The entropy estimate
$\CH_N(p) = \CH_{N,{\bf Y},{\bf Z}}(p)$
introduced by Gy\"{o}rfi and Van Der Meulen \cite{Gyor2}, is then
defined by:
\begin{eqnarray}
\label{entrop}
\CH_N(p)=\frac{1}{[N/2]}\sum_{i=1}^{[N/2]}
\log \hat p_N(Y_i) \I_{\left\{p_N(Y_i)\geq a_N\right\}}
\end{eqnarray}
where $0<a_N<1$ and $\lim_{N\rightarrow \infty} a_N=0.$

\begin{theorem}
      \label{th:main}
Assume that ${\cal H}(f)<\infty$.
For all $n\geq 0$, let
$\BX_N^n$ be a N-sample from $p^n$, the p.d.f.\ of the
MH algorithm at time~$n$, and
consider the kernel density estimate $\widehat{p^n}_N$ given in (\ref{parzen}),
based on the subsample $\BZ^n_N$ defined in (\ref{Z}).  Let
the kernel $K$ be a bounded density, vanishing outside a sphere $S_r$
of radius $r>0$, and set
$h_N=N^{-\alpha}$, $0<\alpha<1/s$.
Consider the entropy estimate $\CH_N$ defined in (\ref{entrop}) with
\begin{eqnarray}
a_N=(\log N)^{-1}.
\label{eq:aN}
\end{eqnarray}
Assume that there are positive constants
$C$, $r_0$, $a$, $A$ and $\epsilon$, such that either:

\begin{itemize}
\item[(i)]  in the case of the Independance Sampler:
$f$, $q$ and $p_0$ satisfy conditions of Proposition~\ref{prop:pnLip};
$q$ satisfies the minoration condition
$q(y)\ge a f(y)$, and~$f$ satisfies the tail condition
\begin{equation}
\label{eq:queuef}
f(y)\leq \frac{C}{||y||^s(\log ||y||)^{2+\epsilon}},
\quad\mbox{for $||y||>r_0$};
\end{equation}

\item[(ii)]  in the general  MH case: $f$, $q$ and $p_0$
satisfy conditions of Proposition~\ref{prop:LipgeneralMH};
$q$ is symmetric ($q(x|y)=q(y|x)$);
$\left|\left| p^0/ f\right|\right|_{\infty} \le A$,
and~$f$ satisfies the tail condition
\begin{equation}
\label{eq:queueGeneral}
f(y) \le {C\over 1+||y||^{s+\epsilon}} .
\end{equation}
\end{itemize}

Then, for all $n\geq 0$,
${\cal H}_N(p^n)\stackrel{a.s.}{\longrightarrow} {\cal H}(p^n)$,
as $N\to\infty$.
\end{theorem}

\begin{proof}
This result uses directly Gy\"orfi and Van Der Meulen's Theorem in
\cite{Gyor2} p.  231.  Conditions (\ref{eq:queuef}) or
(\ref{eq:queueGeneral}) and the
fact that ${\cal H}(f)<\infty$
implies, for all $n\geq 0$, the same conditions on the densities
$p^n$ in either cases (i) or (ii).
Actually, ${\cal H}(f)<\infty$ is a direct consequence of (\ref{taux}) and
of the positivity of $\cal K$.
For the tail condition (\ref{eq:queuef}), case (i), it suffices to notice that
from (\ref{eq:vn}) we have for all $x\in \Omega$:
\begin{eqnarray}
0\leq p^n(x)&\leq& f(x)+\kappa \rho^n f(x)\nonumber\\
&\leq& \frac{C(1+\kappa \rho^n)}{||x||^s(\log ||x||)^{2+\epsilon}}.\nonumber
\end {eqnarray}
The tail condition for the general case (ii) comes directly from the recursive
formula~(\ref{eq:pngeneral}) since
\begin{eqnarray*}
p^1(y) &=& J_0(y)+p^0(y)(1-J(y))
\le \int p^0(x)q(y|x)\alpha(x,y)\,dx + p^0(y)\\
&\le& \int p^0(x)q(y|x){f(y)\over f(x)}\,dx + p^0(y)\\
&\le& A f(y) \int q(x|y)\,dx   + p^0(y) \le 2A f(y).
\end{eqnarray*}
Applying this recursively gives
$$
p^n(y) \le 2^n A f(y) \le {2^n A C \over 1+||y||^{s+\epsilon}},
$$
which is stricter than Gy\"orfi and Van Der Meulen's tail condition.
As to smoothness,
the conditions of our Proposition~\ref{prop:pnLip} for case (i),
and Proposition~\ref{prop:LipgeneralMH} for case (ii) give the
Lipschitz condition of Ivanov and Rozhkova~\cite{Iva} for $p^n$, which
in turn is stricter than Gy\"orfi and Van Der Meulen's smoothness
condition, as stated in Gy\"orfi and Van Der Meulen~\cite{Gyor2}.
\end{proof}

\section{Examples}
\label{sec:simul}

We give in this section several examples for synthetic models,
with target densities which are one and two-dimensional mixtures of gaussian
distributions. The advantage of taking a mixture is that it is an easy
way to build multimodal target densities with ``almost disconnected''
modes, i.e.\ separated modes with regions of low probability in between
(see figures~\ref{fig:truepdf} and~\ref{fig:2dimpdf}).

The difficulty for classical RWMH algorithm is then to properly calibrate
the variance of the random walk to propose jumps under all the modes in a
reasonable amount of time. The difficulty for the IS is to
use a good proposal density~$q$, hopefully allowing sufficient
mass over the modal regions.

It is important to understand that in all the examples below,  the
target density is completely known so that, instead of estimating the
difference $\CK(p_{1}^n,f) - \CK(p_{2}^n,f)$ between any two given
strategies, we are able to estimate directly
$\CK(p_i^n,f)$ for each strategy~$i$ leading to the successive densities
$p_i^n$
separately. Actually, we can compute the strongly consistent estimate
$$
\CK_N(p_i^n,f) = \CH_N(p_i^n)
- {1\over N}\sum_{j=1}^N \log f(X_j^{i,n}) ,
$$
where the $X_j^{i,n}$'s, $j=1,\ldots,N$ are i.i.d.\ $\sim p_i^n$.

We give the results in terms of these estimates, since they provide
easier comparisons and illustrate more clearly the behaviour of our
method.  However, one should keep in mind that in real-size
situations, only the plots of the differences are accessible to
computation. This is not a flaw in the method since clearly,
the better algorithm can be deducted from these plots. For complete
illustration, however,
we have also provided for the first example the plots of
$n\mapsto D(p_1^n,p_2^n,f)$ for comparing three strategies.
The only information {\it not} provided by the plot of the difference is the
``convergence time'' of each chain (in the sense of the convergence
assessment of MCMC, see, e.g., Gilks {\it et al.}~\cite{GRS}). Indeed,
even if the difference goes about zero at time~$n$,
there is always a possibility
that both MH algorithms fail to converge at that time,
with $\CK(p_{1}^n,f) \approx \CK(p_{2}^n,f)$.

Since the models are quite simple here, we could ran a large number of
i.i.d.\ Markov chains to obtain precise estimates.  So we tried up to
$N=1000$ chains for the one-dimensional model, and up to $N=200$ chains
for
the two-dimensional model.  This number of parallel chains can be
reduced without compromising the decision for real-size applications
(we tried $N=100$ chains with satisfactory results).  Note also that
the computing time needed, even with large $N$, is not long since
the duration $n_0$ of the parallel simulation is itself short: the
best algorithm is quickly detected, as shown in the figures below.
Finally, for computing the estimate of the
entropy~(\ref{entrop}), we use a treshold
$a_N = \CO(\log(N)^{-1})$ instead of~(\ref{eq:aN}), to avoid rejection
of too many observations for small values of~$N$.

\subsection{A one-dimensional example}
\label{ss:exsimu}

\begin{figure}[h]
\centering
\includegraphics[height=3truecm, width=8truecm]{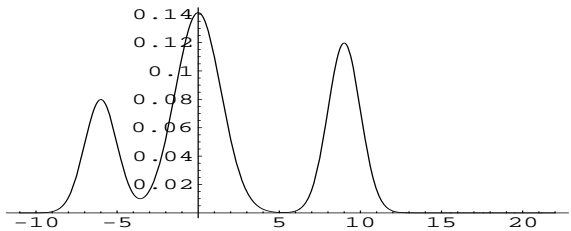}
\caption{\small\sl True mixture density~(\ref{eq:truepdf}).}
\label{fig:truepdf}
\end{figure}

To illustrate the relevance of our comparison of the convergence
rates, we first choose a very simple but meaningful situation,
consisting in MH algorithms for simulation from a mixture of 3
gaussian distributions, with density
\begin{equation}
\label{eq:truepdf}
f(x) = \sum_{i=1}^3 \alpha_{i}\, \varphi(x ; \mu_{i},\sigma^2_{i}),
\end{equation}
where $\varphi(\,\cdot\, ; \mu,\sigma^2)$ is the pdf of
${\cal N}(\mu,\sigma^2)$. The chosen parameters
$\alpha_1=0.5$, $\alpha_2=0.3$, $\mu_1=0$,
$\mu_2=9$, $\mu_3=-6$, $\sigma^2_1=2$, and
$\sigma^2_2=\sigma^2_3=1$,
result in the trimodal pdf
depicted in figure~\ref{fig:truepdf}.

\paragraph{Independence Sampler}

\begin{figure}[t]
\centering
\includegraphics[height=3truecm, width=50truemm]{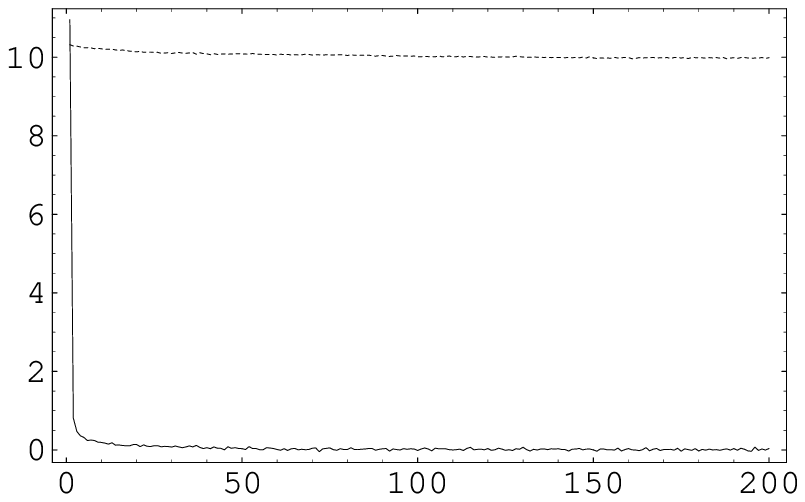}
\includegraphics[height=3truecm, width=50truemm]{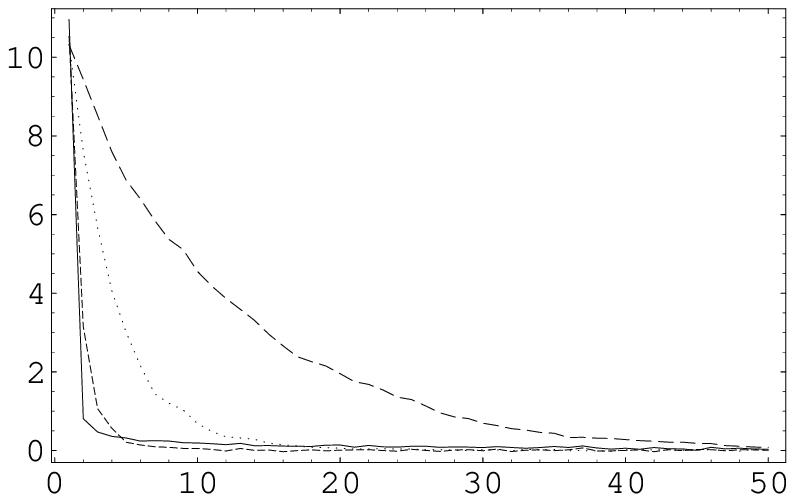}
\caption{\small\sl Plots of
$n\mapsto \CK_N(p^n,f)$ for the IS with a gaussian
proposal $\CN(0,\sigma^2)$.
Left: $\sigma=3$ (solid) vs. $\sigma=1$ (dashed);
Right: $\sigma=3$ (solid) vs. $\sigma=10$ (dashed),
$\sigma=30$ (dotted), $\sigma=100$ (long dashed).}
\label{fig:exIS}
\end{figure}

We first ran the independence sampler with a gaussian proposal
density~$q = \CN(0,\sigma^2)$, for several settings of the variance
parameter.
None of these MH are optimal, since $q$ does not have
modes at $-6$ and $9$, whereas $f$ has. If the variance is too small
(e.g., $\sigma\le1$), the algorithm can ``almost never'' propose jumps
outside, say, $[-4;+4]$, so that it can (almost) never visit the
right  or left modes. The algorithm requires then a
dramatically long time to converge. Our Kullback divergence estimate
reflect this fact (figure~\ref{fig:exIS}, left). For more adapted
settings like, e.g., $\sigma=3$, the algorithm converges faster, and
the convergence again deteriorates as $\sigma$ increases above, say,
10, since the proposal density is then overdispersed
(see figure~\ref{fig:exIS}, right).

\begin{figure}[h]
\centering
\includegraphics[height=3truecm, width=50truemm]{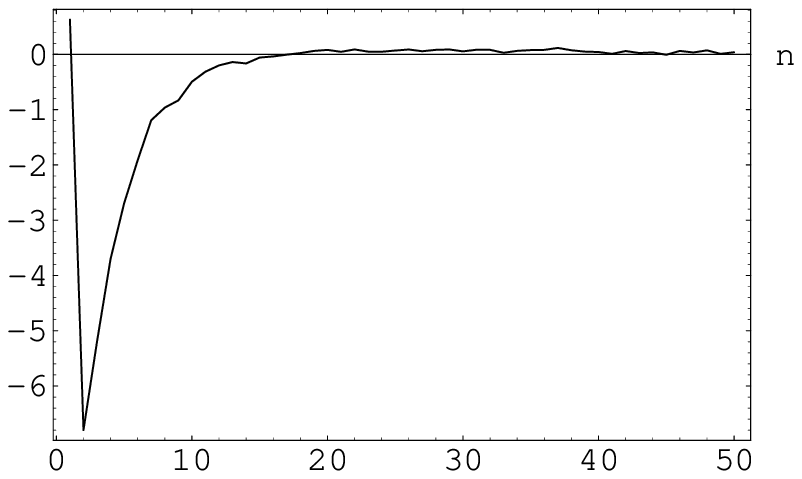}
\includegraphics[height=3truecm, width=50truemm]{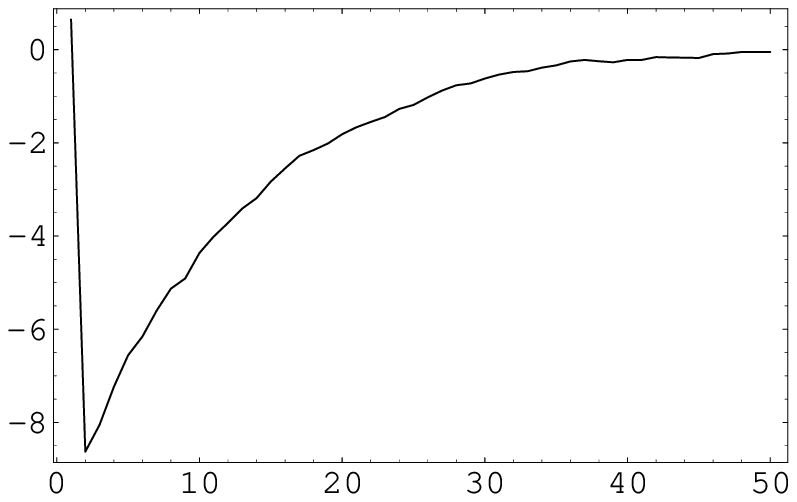}
\caption{\small\sl Plots of the difference
$n\mapsto D(p_1^n,p_2^n,f)$ for the IS with a gaussian
proposal $q_i = \CN(0,\sigma^2_i)$.
Left: $\sigma_1=3$ vs. $\sigma_2=30$;
Right: $\sigma_1=3$ vs. $\sigma_2=100$.}
\label{fig:exISDiff}
\end{figure}

For this example, we also provide
in figure~\ref{fig:exISDiff} two
examples of the plots available in actual situations, i.e.\ that of
$n\mapsto D(p_1^n,p_2^n,f)$. The sign of the plots in both cases, and
for the first iterations,
clearly indicate that, in both cases, the first strategy
($q_1 = \CN(0,3^2)$) is preferable. Morevover, the comparison of the
two plots indicate that $\sigma=100$ is even worse than $\sigma=30$.

\begin{figure}[h]
\centering
\includegraphics[height=3truecm, width=50truemm]{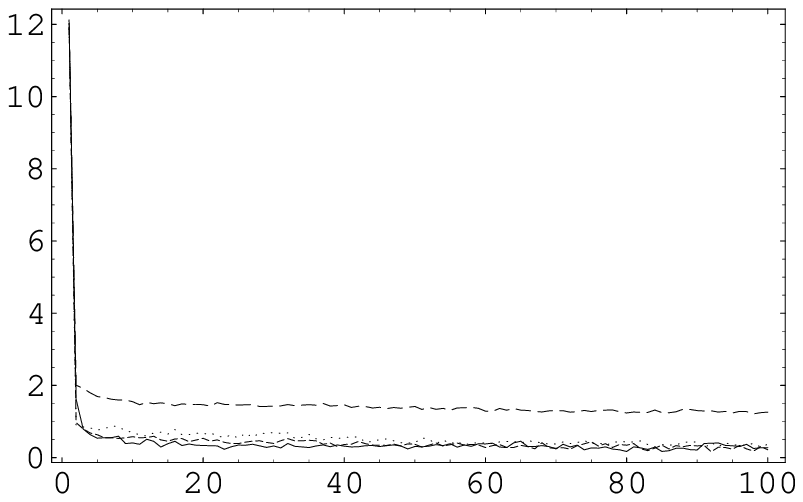}
\includegraphics[height=3truecm, width=50truemm]{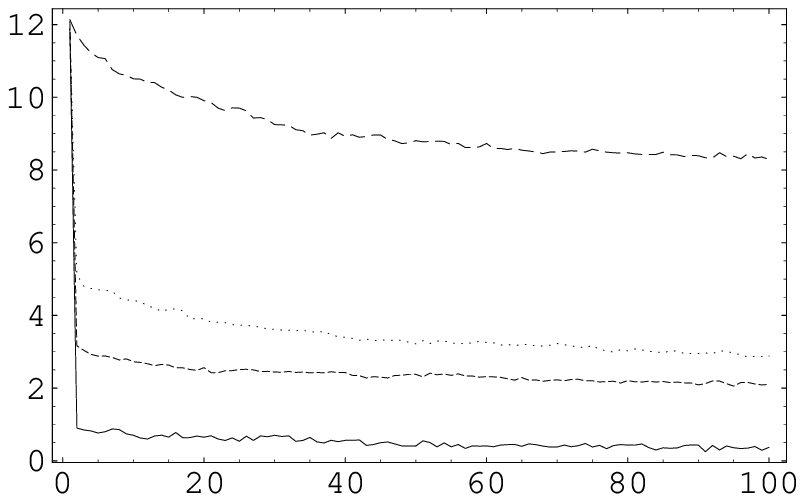}
\caption{\small\sl Plots of
$n\mapsto \CK_N(p^n,f)$ for the IS with a Student
proposal $t(d)$: Left: $d=1$ (solid), $d=2$ (short dashed),
$d=3$ (dotted), $d=10$ (long dashed).
Right: $d=3$ (solid), $d=20$ (short dashed),
$d=50$ (dotted), $d=100$ (long dashed).}
\label{fig:exISstudent}
\end{figure}

To check our estimates with heavy-tailed proposal densities, we also
ran the independence sampler with a Student proposal density
$t(d)$, for $d=1$ (the Cauchy distribution),
up to $d=100$ (for which the Student is almost the normal
distribution). As expected, the algorithms converge faster when they
use Student distributions with heavier tails, since in that case they can
still propose jumps in the left or right mode. When $d\to\infty$, the
proposal converges to the $\CN(0,1)$, and the IS shows the same
behavior as the previous one, with $\sigma=1$ (compare
figure~\ref{fig:exIS}, left with figure~\ref{fig:exISstudent}, right).

\paragraph{RWMH}

\begin{figure}[h]
\centering
\includegraphics[height=3truecm, width=50truemm]{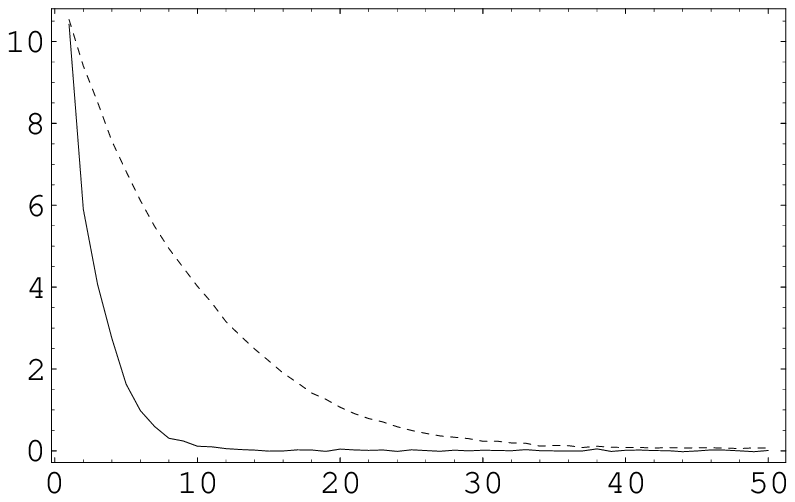}
\includegraphics[height=3truecm, width=50truemm]{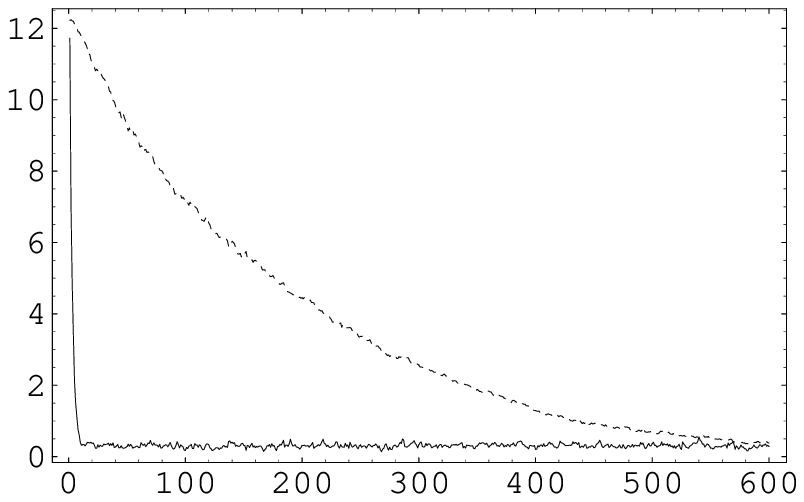}\newline
\caption{\small\sl Plots of
$n\mapsto \CK_N(p^n,f)$ for the RWMH with gaussian proposal
$q(\cdot|x) \equiv \CN(x,\sigma^2)$:
left: $\sigma=1$ (dashed) vs. $\sigma=10$ (solid);
right: $\sigma=0.1$ (dashed) vs. $\sigma=10$ (solid).}
\label{fig:exRW}
\end{figure}

We also ran  on the same example a random walk MH algorithm
with a gaussian proposal $q(\cdot|x) \equiv \CN(x,\sigma^2)$, and
several settings for $\sigma^2$. As expected in view of the region of 
interest for the target~$f$, a good choice is about
$\sigma=10$. For
too small settings (e.g., $\sigma=0.1$), the jumps of the random walk
(of order $\pm 3\sigma$) are
too small, so that the chain needs a
dramatically long time to reach the rightmost mode. This is clearly
indicated by our estimate in figure~\ref{fig:exRW}, right, where up to
$n=600$ iterations are needed for this inefficient algorithm to
converge.

\subsection{A two-dimensional example}
\label{ss:exsimu2d}

\begin{figure}[h]
\centering
\includegraphics[height=4truecm, width=50truemm]{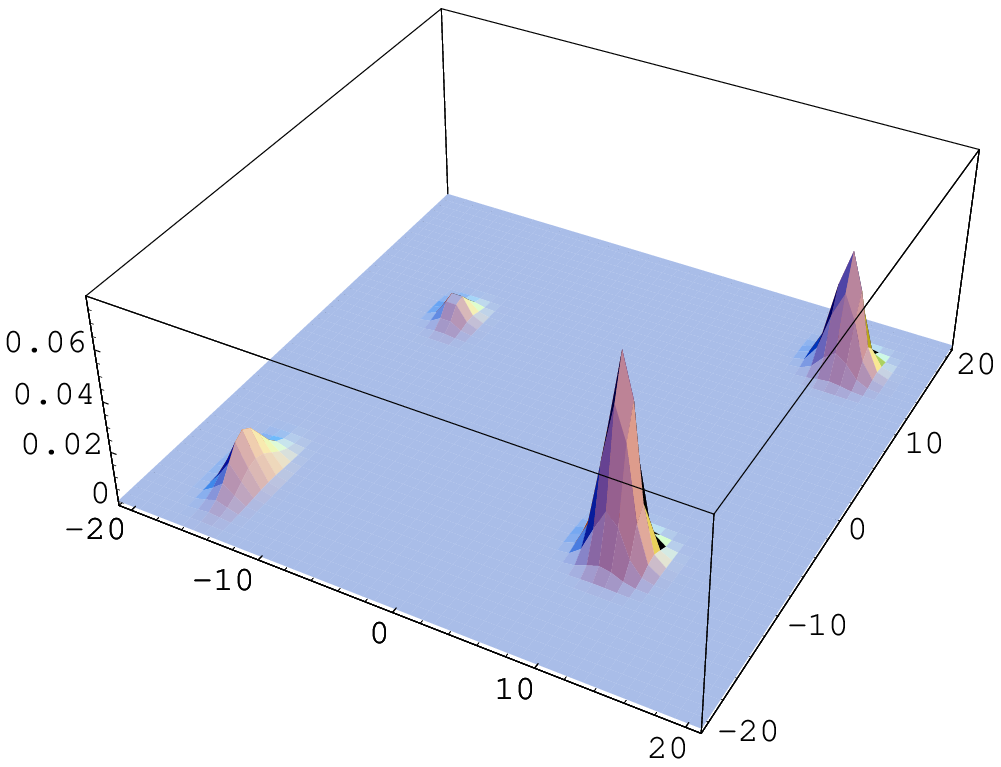}
\includegraphics[height=4truecm, width=50truemm]{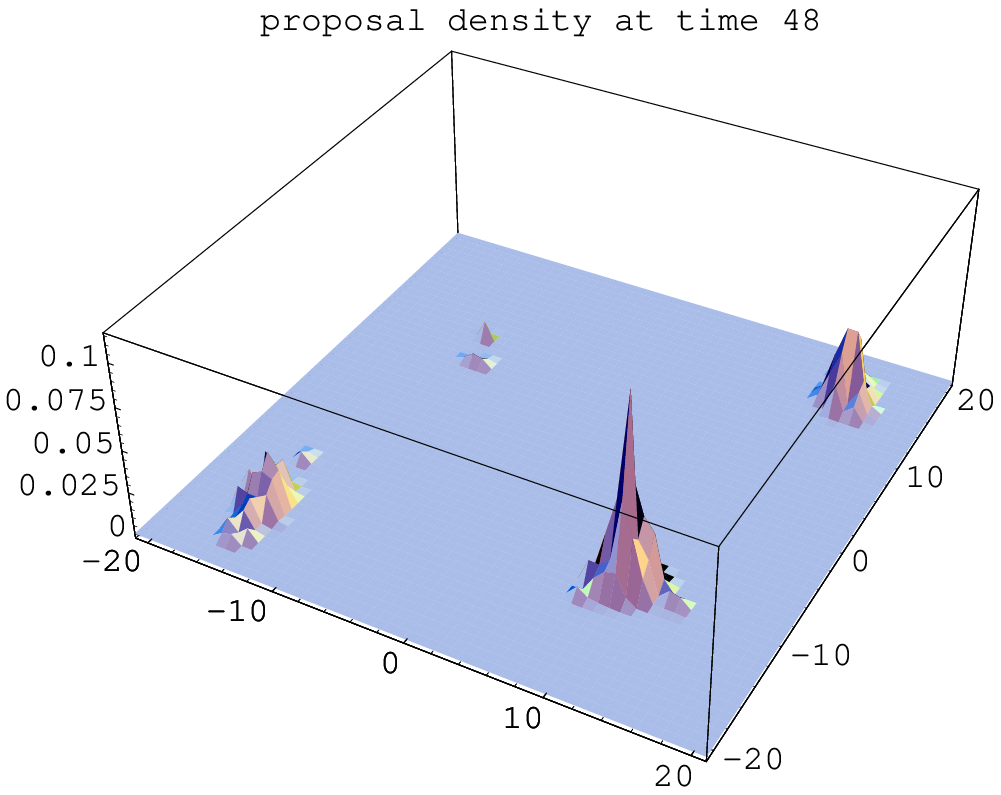}\newline
\caption{\small\sl left: true target pdf;
right: adaptive proposal density.}
\label{fig:2dimpdf}
\end{figure}

We also tried for the target density
a two-dimensional gaussian mixture, depicted in
figure~\ref{fig:2dimpdf}, left
(the true parameters are not given for brevity).
For this example, we compare three ``good'' strategies
of different types: (i) An IS with a uniform proposal density over the
compact $[-20;20]^2$; this algorithm is ``almost geometric'' since the
mass of the tails of~$f$ outside the compact are negligible
(the minoration condition $q\ge a f$ is fulfilled on the compact).
(ii) A RWMH with a bivariate gaussian proposal
$q(\cdot|x) = \CN(x,\sigma^2 I)$, with a ``good'' setting
$\sigma=17$, founded using our Kullback divergence method.
(iii) An adaptive MH algorithm following the ideas in
Chauveau and Vandekerkhove~\cite{CV1}. In short, parallel chains
started with the IS (i) are
ran and kernel density estimates are built at some specified
times using the past of these i.i.d. chains. For example,
the proposal density built at time $n=48$ is depicted in
figure~\ref{fig:2dimpdf}, right.

\begin{figure}[h]
\centering
\includegraphics[height=30truemm, width=50truemm]{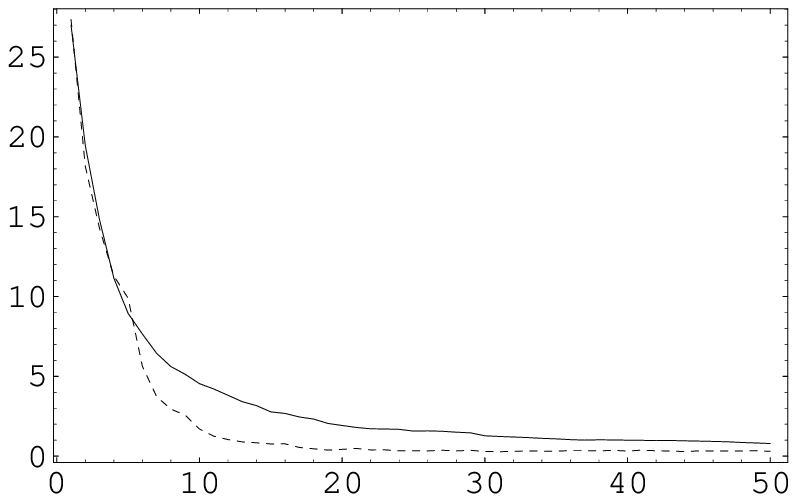}
\includegraphics[height=30truemm, width=50truemm]{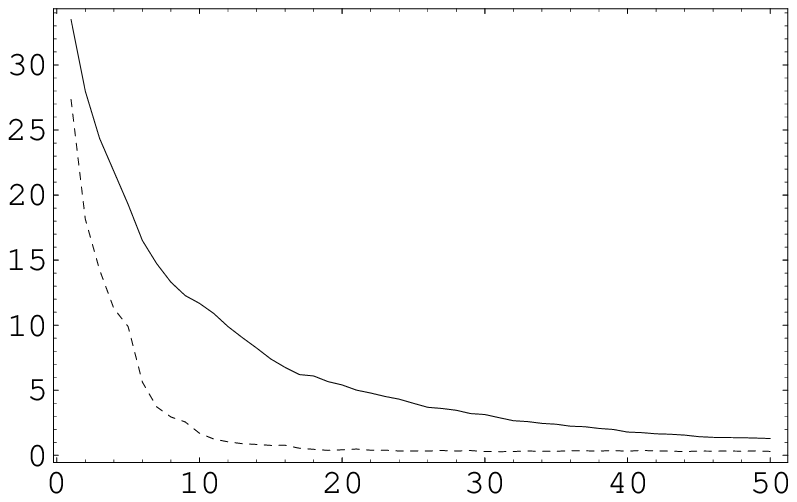}\newline
\caption{\small\sl Plots of
$n\mapsto \CK_N(p^n,f)$ based on $N=200$ chains
for the IS with adaptive proposal (dashed)
vs. left: IS with $q =$ Uniform distribution;
right: RWMH with $q(\cdot|x) = \CN(x,\sigma^2 I)$, $\sigma=17$.}
\label{fig:2dimK}
\end{figure}

The estimates of $\CK(p^n,f)$ for this setup (and $N=200$ chains) are
given in figure~\ref{fig:2dimK}.  As expected, the IS with the uniform
proposal density performs better than the calibrated RWMH. But the
adaptive proposal is even better than the two others.  The times at
which the adaptive proposal density is updated ($n=5$ and $9$), are
even visible in the plot of $\CK(p^n,f)$ for this strategy, which
means that it has an immediate effect on this divergence.

\section{Conclusion}

We have proposed a methodology to precisely quantify and compare
several MCMC simulation strategies only on the basis of the simulations
output. A procedure for applying our
method in practice has been given in Section~\ref{sec:int}.

A novelty is that this methodology is based upon the use of the
relative entropy of the successive densities of the MCMC algorithm.  A
consistent estimate of this entropy in the MH setup has been proposed,
and general conditions insuring its convergence have been detailed.
Indeed, the conditions of propositions~\ref{prop:pnLip}
and~\ref{prop:LipgeneralMH} are difficult to verify in practice.
However, the theoretical study of Section~\ref{sec:smooth} has been
developed to support the idea that, if the ``ingredients'' of the MH
algorithm ($q$, $p^0$ and~$f$) have sufficient tails and smoothness
conditions, then one can reasonably expect that the successive
densities $p^n$, $n\ge1$, of the MH algorithm will also satisfy these
conditions, so that our usage of the estimates of the $\CH(p_i^n)$'s
will indicate the most efficient MH algorithm to use.

Our methodology is an alternative to the adaptive MCMC strategies,
which represent an active field of the current literature on this
field.  The advantage of our approach is that it avoids the
difficulties associated to adaptation, namely the preservation of the
convergence to the desired limiting distribution, and the theoretical
guarantee that the adaptive algorithm will perform better than any
other heuristical approach.

Most importantly, our method can also be used as a global criterion to
compare, on specific cases, the incoming new (eventually adaptive)
MCMC strategies against existing simulation methods.  Note in addition
that, if the comparisons are done on simulated situations where $f$ is
entirely known, our approach gives directly the estimate of $\CK(p^n,f)$
for each strategy instead of the difference between two methods.
We used for this purpose in example~\ref{ss:exsimu2d}, and
the need for such a global comparison criterion on simulated
situations is apparent in recent literature, as e.g.~in 
Haario {\it et. al\/}~\cite{HSTpp} and~\cite{HST}.

\section{Appendix}

The purpose of this appendix is to show that some of the conditions
required in Proposition~\ref{prop:pnLip} and
Proposition~\ref{prop:LipgeneralMH}, which look difficult to check in
actual situations, are satisfied at least in simple situations, for
classically used MH algorithms in the one-dimensional case
i.e.\ when $x\in\bR$.

\subsection{The one-dimensional independence sampler case}
\label{ss:ApIS}

In the IS case, the difficult conditions are
conditions~(\ref{eq:condLipalpha}) and~(\ref{eq:condcq}) of
Lemma~\ref{lem:LipI}.
These conditions are simpler to handle
in the one-dimensional case.
First, note that we can prove under additional conditions
the derivability of $p^n$ for all
$n\ge1$ (that proof is not given since we are not using it here).
In the one-dimensional case, when $q$ and $f$ are in addition
derivable, and have non-oscillating tails,
assumption (A) leads to
\begin{equation}
\label{eq:H}
\exists m_1 < m_2 :
\forall x < m_1, h'(x) < 0, \mbox{ and }
\forall x > m_2, h'(x) > 0.
\end{equation}
For a fixed $x\in\bR$, there exists by~(\ref{eq:H}) a compact set
$K(x) = [a(x) , b(x)]$ such that
\begin{enumerate}
          \item[(i)]  $[m_1,m_2] \subseteq K(x)$;

          \item[(ii)]  $h(a(x)) = h(b(x)) = h(x)$;

          \item[(iii)]  for any $y\notin K(x)$, $h(y) \ge h(x)$.
\end{enumerate}
As in the general case, this entails that
$\forall y \notin K(x)$, $\alpha(y,x) = 1$. If we have the Lipschitz
condition on $K(x)$:
$$
\forall y,z, \quad |\alpha(y,x) - \alpha(z,x)| \le c(x) |y-z|,
$$
the expression of $c(x)$ can be precised
\begin{equation}
          c(x) = \sup_{y\in K(x)} \left|
          {\partial \phi(y,x)\over \partial y}\right| < \infty.
          \label{eq:c}
\end{equation}
and Lemma~\ref{lem:LipI} holds if the integrability
condition~(\ref{eq:condcq}) is satisfied.
Note that $|a(x)|$ and $b(x)$ both go to $+\infty$ as
$|x|\to\infty$; in particular, $b(x)=x$ for $x>m_2$.
Hence $c(x)\to\infty$ as $|x|\to\infty$, and
condition~(\ref{eq:condcq}) is
not always true, but merely depends on the relative decreasing rate of
the tails of~$q$ and~$f$.

For an illustrative example, assume that the tails of~$f$
are of order $x^{-\beta}$, and the
tails of $q$ are of order $x^{-\alpha}$. Satisfying assumption A requires that
$\beta>\alpha$. Now, one can always use the fact that
$$
c(x) \le \sup_{y\in \bR} \left|
          {\partial \phi(y,x)\over \partial y}\right|,
$$
so that if $\beta-1 < \alpha < \beta$, then
$c(x)$ is of order $x^{\alpha-\beta}$ for large $x$
and~(\ref{eq:condcq}) is satisfied. The condition
$\alpha\in[\beta-1,\beta]$ states that the tails of~$q$
should be ``not too heavy'', compared with the tails of~$f$.
This requirement is obviously stronger than
what is needed, but more precise conditions require some analytical
expression of $c(x)$ for $x\notin [m_1,m_2]$, and this expression
depends on $a(x)$ and $h'$.

Fortunately, condition~(\ref{eq:condcq}) is satisfied in much more
general settings. For instance, consider situations where~$f$
and~$q$ are both symmetric w.r.t.\ $0$, so that $K(x) = [-|x|,|x|]$
for $x$ outside $[m_1,m_2]$,
and $c(x)$ can be expressed in closed form.
Then it is easy to verify that~(\ref{eq:condcq}) holds for, e.g.,
$f \equiv \CN(0,1)$ and $q\equiv t(d)$, the Student distribution with
$d$ degrees of freedom, for $d\ge 2$ (even if, for $d=2$ the tails of
$q$ are of order $x^{-3}$).
In this example, the proposal density has tails
much more heavier than~$f$, but Lemma~\ref{lem:LipI} holds
i.e., $I$ is still Lipschitz.

\subsection{The one-dimensional general MH case}
\label{ss:ApGeneral}

In the general MH case, the difficult conditions are conditions~(ii)
and~(iii) of Proposition~\ref{prop:LipgeneralMH}. Our aim is to show
that these conditions hold in the simple RWMH case with gaussian
proposal density. In order to obtain a tractable case, let
$q(y|x)$ be the p.d.f.\ of the gaussian
$\CN(x,1)$, and $f$ be the density of the target distribution $\CN(0,1)$.

For condition~(ii) we have to prove that
$q(\cdot|x)\alpha(x,\cdot)$ is $c(x)$-Lipschitz, with
$\int p^n(x) c(x)\,dx < \infty$. Here $q(y|x)=q(x|y)$, so that
$$
\alpha(x,y) = 1 \wedge {f(y)\over f(x)} \le {f(y)\over f(x)},
$$
which is a truncated function such that, for any $x$,
$\lim_{|y|\to\infty}\alpha(x,y) = 0$. In other words,
both $\alpha(x,y)$ and $q(y|x)\alpha(x,y)$ have tails behavior
for large~$y$.
The non-troncated function
$\varphi_x(y) = q(y|x)f(y)/f(x)$ is then Lipschitz, with
$$
c(x) = \sup_{y\in\bR}\left| \varphi_x'(y) \right|.
$$
A direct calculation (feasible in this simple case) gives
$c(x) \propto \exp(x^2-2)/4$.
Since to ensure the tails conditions of the successive densities
$p^n$ we have to assume that the initial distribution itself has tails lighter
or equal to that of~$f$ (i.e.\ that $||p^0/f||_{\infty} < A$, see
Theorem~\ref{th:main}) then by
the recursive
definition of $p^n$ we have, as in the proof of Theorem~\ref{th:main},
$p^n(y) \le 2^n A f(y)$, so that
$\int p^n(x) c(x)\,dx < \infty$, i.e.\ condition~(ii) of
Proposition~\ref{prop:LipgeneralMH} holds.

We turn now to condition~(iii) of Proposition~\ref{prop:LipgeneralMH},
i.e.\ we have to show that $J(y)$ given by~(\ref{eq:J}) is Lipschitz.
For fixed $y,z \in\bR$,
$$
|J(y)-J(z)| \le
\int |q(x|y)\alpha(y,x) - q(x|z)\alpha(z,x)|\,dx.
$$
As for the IS case, we need a precise study of the truncated function
here. We assume first that $z>y>0$. Since $q$ is symmetric,
$$
\alpha(y,x) = {f(x)\over f(y)}\wedge 1,
$$
and we can define two compact sets $K(y)$ and $K(z)$ by
$$
K(t) = \{x\in\bR:\alpha(t,x)=1\} = \{x \in\bR: f(x)\ge f(t)\}
$$
which, in the present situation, are just
$K(y) = [-y,y]$, $K(z)=[-z,z]$, and satisfy $K(y)\subset K(z)$.
Hence
\begin{eqnarray*}
|J(y)-J(z)| &\le& \int_{K(y)} |q(x|y) - q(x|z)|\,dx \\
& &+\int_{K(z)\setminus K(y)}
\left|q(x|y){f(x)\over f(y)} - q(x|z)\right|\,dx \\
& &+\int_{K(z)^c}
\left|q(x|y){f(x)\over f(y)} - q(x|z){f(x)\over f(z)}\right|\,dx,
\end{eqnarray*}
where $K(z)^c = \bR\setminus K(z)$.
Using the mean value theorem, the first term can be written
\begin{eqnarray}
\int_{K(y)}
|q(x|y) - q(x|z)|\,dx &\le& \int |q(x|y) - q(x|z)|\,dx \nonumber\\
&\le& |y-z|
\int {|x-y^*|\over 2\pi}\exp\left(-(x-y^*)^2/2\right) \,dx \nonumber\\
&\le&  \sqrt{2\over\pi}\,|y-z|, \label{eq:term1}
\end{eqnarray}
where the last inequality comes from the absolute first moment of the
normal density.

For the second term, consider first the integral on the right side
of $K(z)\setminus K(y)$, that is
$\int_y^z \left|\varphi_{y,z}(x)\right|\,dx$, where
$$
\varphi_{y,z}(x) = q(x|y){f(x)\over f(y)} - q(x|z).
$$
In this simple setting, it is easy to check that $\varphi_{y,z}(\cdot)$
is a bounded function, monotonically decreasing from
$\varphi_{y,z}(y) = \delta - q(y|z) >0$ to
$\varphi_{y,z}(z) = q(y|z)- \delta <0$, where
$\delta = q(y|y)$ is
the value of the gaussian density at its mode. Hence
\begin{equation}
      \label{eq:term2}
\int_y^z \left|q(x|y){f(x)\over f(y)} - q(x|z)\right|\,dx
\le \delta |y-z|.
\end{equation}
The symmetric term $\int_{-z}^{-y} \left|\varphi_{y,z}(x)\right|\,dx$
is handled in a similar way.

The third term can in turn be decomposed into
\begin{eqnarray*}
\int_{K(z)^c}
\left|q(x|y){f(x)\over f(y)} - q(x|z){f(x)\over f(z)}\right|\,dx
&\le&
Q\int_{K(z)^c}
\left|{f(x)\over f(y)} - {f(x)\over f(z)}\right|\,dx \\
&+& \int_{K(z)^c}
\left|q(x|y) - q(x|z)\right|\,dx ,
\end{eqnarray*}
where, as in Proposition~\ref{prop:LipgeneralMH}, $Q=||q||_{\infty}$,
and since $\sup_{x\in K(z)^c} |f(x)/f(z)| = 1$.
Using the mean value theorem as for the first term,
\begin{equation}
      \label{eq:term31}
\int_{K(z)^c}
\left|q(x|y) - q(x|z)\right|\,dx \le \sqrt{2\over\pi} \, |y-z|.
\end{equation}
Finally,
\begin{eqnarray}
\int_{K(z)^c}
\left|{f(x)\over f(y)} - {f(x)\over f(z)}\right|\,dx
&=& 2 \int_z^\infty
\left|{f(x)\over f(y)} - {f(x)\over f(z)}\right|\,dx \nonumber\\
&\le&
2\left|{1\over f(y)} - {1\over f(z)}\right|
\int_z^\infty f(x)\,dx \nonumber\\
&\le& 2\sqrt{2\pi} z e^{z^2/2}
{e^{-z^2}\over z+\sqrt{z^2+4/\pi}}|y-z|, \label{eq:Dz}\\
&\le& D |y-z|, \label{eq:term32}
\end{eqnarray}
where the left term in~(\ref{eq:Dz}) comes from the mean value theorem applied
to the function $1/f(\cdot)$, the rightmost term in~(\ref{eq:Dz})
is a well-known bound of the tail of the normal distribution, and
$$
D = \sup_{z\in\bR}\left|
2\sqrt{2\pi} z e^{z^2/2}
{e^{-z^2}\over z+\sqrt{z^2+4/\pi}}
\right| < \infty.
$$
Collecting (\ref{eq:term1}), (\ref{eq:term2}), (\ref{eq:term31}) and
(\ref{eq:term32}) together shows that
$$
|J(y)-J(z)| \le k |y-z| \quad\mbox{for $z>y>0$ and $0<k<\infty$}.
$$
The other cases are done similarly, so
that $J(\cdot)$ is Lipschitz.


\bigskip
\noindent 
{\small{\bf Corresponding author}\\
Didier Chauveau\\
Laboratoire MAPMO - UMR 6628 - F\'ed\'eration Denis Poisson\\
Universit\'e d'Orl\'eans\\
BP 6759, 45067 Orl\'eans cedex 2, FRANCE.\\
Email: {\tt didier.chauveau@univ-orleans.fr}
}


\begin{thebibliography}{99}

\bibitem{Ahm1}
Ahmad, I. A. and Lin, P. E. (1976),
A nonparametric estimation of the entropy for absolutely continuous
distributions,
{\it IEEE Trans. Inform. Theory}, vol. 22, 372--375..

\bibitem{Ahm2}
Ahmad, I. A. and Lin, P. E. (1989),
A nonparametric estimation of the entropy for absolutely continuous
distributions,"
{\it IEEE Trans. Inform. Theory}, vol. 36, 688--692.


\bibitem{AP05}
Atchad\'e, Y.F., and Perron, F. (2005),
Improving on the independent Metropolis-Hastings algorithm,
{Statistica Sinica}, {\bf 15}, no 1,3--18.


\bibitem{AR05}
Atchad\'e, Y.F., and Rosenthal, J. (2005),
On adaptive Markov chain Monte Carlo algorithms,
{\it Bernoulli},{\bf 11}(5), 815--828.

\bibitem{BBN}
Ball, K., Barthe, F. and Naor, A., (2003),
Entropy jumps in the presence of a spectral gap,
{\it Duke Mathematical Journal}, {\bf 119}, 1, 41--63.


\bibitem{Bil}
Billingsley (1995),
{\it Probability and Measure}, 3rd Edition, Wiley, New York.


\bibitem{CV1}
Chauveau, D. and Vandekerkhove, P. (2002),
Improving convergence of the Hastings-Metropolis algorithm
with an adaptive proposal,
{\it Scandinavian Journal of Statistics}, {\bf 29}, 1, 13--29.


\bibitem{CV2}
Chauveau, D. and Vandekerkhove, P. (2004),
A Monte Carlo estimation of the entropy for Markov chains,
{\it preprint}.



\bibitem{Delmo1}
Del Moral P., Ledoux M., Miclo, L., (2003),
Contraction properties of Markov kernels.
{\it Probab. Theory and Related Fields},
{\bf 126}, pp. 395--420.


\bibitem{Dev83}
Devroye, L. (1983),
The equivalence of weak, strong and complete convergence in $L^1$ for 
kernel density estimates,
{\it Ann. Statist.}, {\bf 11}, 896--904.


\bibitem{Dmitri1}
Dmitriev, Y. G., and Tarasenko, F. P. (1973),
On the estimation of
functionals of the probability density and its derivatives,
{\it Theory Probab, Appl.} {\bf 18}, 628--633.

\bibitem{Dmitri2}
Dmitriev, Y. G., and Tarasenko, F. P. (1973),
On a class of
non-parametric estimates of non-linear functionals of density,
{\it Theory Probab, Appl.} {\bf 19}, 390--394.


\bibitem{DGMR05}
Douc, R., Guillin, A., Marin, J.M. and Robert, C.P. (2006)
Convergence of adaptive sampling schemes,
{\it Ann. Statist.}, to appear.

\bibitem{Dude}
Dudevicz, E. J. and Van Der Meulen, E. C. (1981),
Entropy-based tests
of uniformity, {\it J. Amer.  Statist.  Assoc.}, {\bf 76} 967--974.

\bibitem{Egger}
Eggermont, P. P. B. and LaRiccia, V. N. (1999),
Best asymptotic Normality of the Kernel Density Entropy Estimator for
Smooth Densities,
{\it IEEE trans. Inform. Theory,} vol. 45, no. 4, 1321--1326.

\bibitem{GasJor}
G\aa semyr, J. (2003),
On an adaptive version of the Metropolis-Hastings algorithm with
independent proposal distribution,
{\it Scand. J. Statist.}, {\bf 30}, no.\ 1, 159--173.



\bibitem{GelSmi}
Gelfand, A.E. and Smith, A.F.M. (1990),
Sampling based approaches to calculating marginal densities.
{\it Journal of the American Statistical Association}
{\bf 85}, 398--409.


\bibitem{GelSah}
Gelfand, A.E. and Sahu, S.K. (1994),
On Markov chain Monte Carlo acceleration,
{\it Journal of Computational and Graphical Statistics}
{\bf 3}, 261--276.


\bibitem{GemGem}
Geman, S. and Geman, D. (1984),
Stochastic relaxation, Gibbs distributions and
the Bayesian restoration of images.
{\it IEEE Trans. Pattern Anal. Mach. Intell.} {\bf 6}, 721--741.


\bibitem{GRG}
Gilks, W.R., Roberts, G.O. and George, E.I. (1994),
Adaptive direction sampling,
{\it The statistician}, {\bf 43}, 179--189.


\bibitem{GRS}
Gilks, W.R., Richardson, S. and Spiegelhalter, D.J. (1996),
{\it Markov Chain Monte Carlo in practice}.
Chapman \& Hall, London.


\bibitem{GiRoSa}
Gilks, W.R., Roberts, G.O. and Sahu, S.K. (1998),
Adaptive Markov chain Monte carlo through regeneration,
{\it Journal of the American Statistical Association}
{\bf 93}, 1045--1054.


\bibitem{Gyor1} 
Gy\"{o}rfi, L. and Van Der Meulen, E. C. (1987),
Density-free convergence properties of various estimators of the entropy,
{\it Comput. Statist. Data Anal.}, 5, 425--436.

\bibitem{Gyor2}
Gy\"{o}rfi, L. and Van Der Meulen, E. C. (1989),
An entropy estimate based on a kernel density estimation,
{\it Colloquia Mathematica societatis J\'anos Bolyai 57.
Limit Theorems in Probability and Statistics P\'ecs (Hungary)},
229--240.




\bibitem{HSTpp}
Haario, H., Saksman, E and Tamminen, J. (1998),
An adaptive Metropolis Algorithm,
{\it Report}, Dpt. of mathematics, University of Helsinki,
Preprint.

\bibitem{HST}
Haario, H., Saksman, E and Tamminen, J. (2001),
An adaptive Metropolis Algorithm,
{\it Bernouilli} {\bf 7}, 2, 223--242.


\bibitem{Has}
Hastings, W.K. (1970),
Monte Carlo sampling methods using Markov Chains and their
applications,
{\it Biometrika} {\bf 57}, 97--109.

\bibitem{Hol}
Holden, L. (1998),
Geometric Convergence of the Metropolis-Hastings Simulation
Algorithm,
{\it Statistics and Probabilitiy Letters}, {\bf 39}, 1998.


\bibitem{Iva}
Ivanov, A. V. and Rozhkova, M.N. (1981),
Properties of the statistical estimate of the entropy of a random
vector with a probability density (in Russian),
{\it Probl. Peredachi Inform}, {\bf 17}, 33-43. Translated into
English in
{\it Problems Inform. Transmission}, {\bf 17}, 171--178.


\bibitem{MenTwe} Mengersen, K.L. and Tweedie, R.L. (1996),
Rates of convergence of the Hastings and Metropolis algorithms.
{\it Ann. Statist.} {\bf 24}, 101--121.


\bibitem{Miclo}
Miclo, L. (1997)
Remarques sur l'hypercontractivit\'e et l'\'evolution de l'entropie
des cha\^{\i}nes de Markov finies.
{\it S\'eminaire de Probabilit\'es XXXI, Lecture Notes in Mathematics},
Springer, 136-168.

\bibitem{Mi01}
Mira, A. (2001),
Ordering and improving the performance of Monte Carlo Markov chains,
{\it Statistical Science}, {\bf 16}, 340--350.


\bibitem{Mokka}
Mokkadem, A. (1989),
Estimation of the entropy and information of
absolutely continuous random variables,
{\it IEEE Trans.  Inform. Theory} {\bf 23} 95--101.

\bibitem{PG05}
Pasarica, C., and Gelman, A. (2005),
Adaptively scaling the Metropolis algorithm using squared jumped
distance,
{\it Technical Report, Columbia University, New York}. 


\bibitem{Ri06}
Rigat, F. (2006),
Markov chain Monte Carlo inference using parallel hierarchical
sampling,
{\it Technical Report, Eurandom, Netherland}.

\bibitem{RoRo}
Roberts, G.O. and Rosenthal, J.S. (2001),
Optimal scaling for various Metropolis-Hastings algorithms,
{\it Statistical Science}, {\bf 16}, 351--367.




\bibitem{RobTwe}
Roberts, G.O. and Tweedie, R.L. (1996),
Geometric convergence and Central Limit Theorems
for multidimensional Hastings and Metropolis algorithms.
{\it Biometrika} {\bf 83}, 95--110.

\bibitem{TW87}
Tanner, M. and Wong, W. (1987),
The calculation of posterior distributions by data augmentation.
{\it J. Am. Stat. Assoc.}, {\bf 82}, 528--550.


\bibitem{Taras}
Tarasenko, F. P. (1968),
On the evaluation of an unknown probability
density function, the direct estimation of the entropy from
independent observations of a continuous random variable, and the
distribution-free entropy test of goodness-of-fit,
{\it Proc.  IEEE.}, {\bf 56} 2052--2053.


\end{thebibliography}
\end{document}